\documentclass[10pt]{amsart}

\newtheorem{theorem}{Theorem}[section]

\theoremstyle{definition}

\theoremstyle{remark}

\numberwithin{equation}{section}

\def\pa{\partial}

\def\na{\nabla}
\def\ve{\varepsilon}

\newcommand\del{\delta}

\begin{document}

\title[Obstruction Function]{An Obstruction Function to the Existence of Complex Structures}

\author[Jun Ling]{Jun \underline{LING}}
\address{Department of Mathematics, Utah Valley University, Orem, Utah 84058}
\email{lingju@uvu.edu}



\date{Dec 14, 2018}



\begin{abstract}
	We construct a function  for  almost-complex manifolds.
	Non-vanishing of the function for the almost-complex structure implies the almost-complex structure  is not integrable.
	Therefore the constructed function is an obstruction for the existence of complex
	structures from the almost-complex structure. It is a function, instead of a tensor.
\end{abstract}

\maketitle

\section{Introduction}\label{sec-intro}\label{sec-intro}

Given  $M^n$  a smooth manifold, an  interesting question is that does there exist any complex manifold structure on $M^n$ that makes $M^n$ a complex manifold? 
For a given complex manifold, the underline complex structure  gives a canonical almost-complex structure. In the study of existence or nonexistence of complex manifold structure for a manifold,  one naturally looks at existence of almost-complex structure first, and  should any exists,  and check  whether or not  it can be "integrated" to a complex structure. 

An almost-complex structure  $J$  on smooth manifold $M^n$  is an endomorphism of the tangent
bundle $TM$  with $J^2=-1$. It is known that if $M$ has an almost-complex structure, then $M$ has even dimension $n$ and $M$ is orientable.  Nijenhuis tensor $N_J$ for the almost complex structure $J$ is given by the following equation
\begin{equation}\label{Nijien}
N_J(X,Y)=[JX,JY]-J[X,JY]-J[JX,Y]-[X,Y],
\end{equation}
for all smooth vector fields $X,Y$.
The celebrated Newlander-Nirenberg theorem  \cite{NN} 
implies that $N_J=0$ if and only if $J$ 
is a canonical almost-complex structure of a complex manifold. 
An almost-complex structure $J$ is called integrable if $N_J$ vanishes.
So in studying  existence or nonexistence of complex manifold structure on a manifold, one often studies  existence or nonexistence of integrable almost-complex structures, or equivalently studies almost-complex structures and the related Nijenhuis tensor
vanishing or non-vanishing property. 

On the other hand, though it is great that vanishing or non-vanishing  Nijenhuis tensor determines whether the almost-complex structure is integrable or not,  it is hard to check the vanishing status for the manifolds other than spheres, due the nature of tensor. 
It would be nice to give some sufficient condition for existence or to give some obstruction for non-existence, of complex structures. In this paper we give an  obstruction that is a function $L_J$, instead of a tensor. Non-vanishing function $L_J$ implies non-vanishing Nijenhuis tensor  $N_J$. To my knowledge, this is the first of obstruction ever appeared. That opens a door for the research of non-existence of complex structures. 
We construct the obstruction function in  Section \ref{sec-LJ}.

\section{An function obstruction for the existence of complex structures}\label{sec-LJ}

We take the convention that we sum on duplicated index in this paper, unless otherwise stated. 
Our first result is

\begin{theorem}\label{thm1}
	If $J$ is an almost-complex structure on a smooth manifold,
	and if
	\begin{equation}\label{LJ}
	L_J:=-\pa_j(J^i_lJ_l^k)\pa_i{J_k^j}.
	\end{equation}
	is not zero at some point of $M$,
	Then $J$ is not integrable.
\end{theorem}
\proof
First take a metric $g$ on the smooth manifold $M$ that always exists. ($L_J$ is independent of $g$ and is in terms of almost-complex structure itself in simple form).
We define a $(4,0)$ tensor $N(X,Z,Y,W)$ by 
\[
N(X,Z,Y,W)
\]
\[:=\frac14\Big\{\left\langle JN_J\Big(N_J(X,Z),Y\Big),W\right\rangle_g
+\left\langle JN_J\Big(N_J(Y,Z),X\Big),W\right\rangle_g
\]
\[
+\left\langle JN_J\Big(N_J(X,W),Y\Big),Z\right\rangle_g
+\left\langle JN_J\Big(N_J(Y,W),X\Big),Z\right\rangle_g\Big\}.
\]
where $N_J$ is the one in (\ref{Nijien}) and  $\langle\cdot,\cdot\cdot \rangle_g =g(\cdot,\cdot\cdot)$.

Take the trace of the  first argument and the third argument and then take the trace
of the second and the fourth argument at a point to get a number for that point
\[
g^{ij}g^{kl}N(\ve_i,\ve_k,\ve_i,\ve_k),
\] 
where $\{\ve_{i}=\frac{\pa}{\pa x^i}\}_{i=1}^n$  is a local frame of local coordinate system  $\{x^i{i}\}_{i=1}^n$, $g_{ij}=g(\frac{\pa}{\pa x^i}, \frac{\pa}{\pa x^j})$, matrix $(g^{ij})=(g_{ij})^{-1}$. 
We write $\pa_i=\frac{\pa}{\pa x^i}$, $\na_i=\na_{\frac{\pa}{\pa x^i}}$
$\langle\cdot,\cdot\cdot \rangle_g =g(\cdot,\cdot\cdot)$ for convenience.

If we can show
\begin{equation}\label{LJ1}
g^{ij}g^{kl}N(\ve_i,\ve_k,\ve_i,\ve_k)=L_J
\end{equation}
then $L_J$ is constructed from Nijenhuis tensor  $N_J$ by contractions. Non-vanishing function $L_J$ implies
non-vanishing tensor $N_J$ since $N_J=0$ implies $L_J=0$. Therefore the theorem follows.

We now prove (\ref{LJ1}).

We calculate at a point with normal coordinates. So at the point $g_{ij}=\del_{ij}$,
$\del_{ii}=1$, $\del_{ij}=0$ if $i\not=j$.

Note that $N(X,Z,X,Z)=\left\langle JN_J\Big(N_J(X,Z),X\Big),Z\right\rangle_g$.
therefore
\[
L_j:=N(e_i,e_k,e_i,e_k)=\left\langle JN\Big(N(e_i,e_k),e_i\Big),e_k\right\rangle_g
=N_{ik}^rN_{ri}^sJ_s^k,
\]
where $N_{ij}^k$ are given by $N_J(\pa_i,\pa_j)=N_{ij}^k\pa_k$. 
It is easy to see
\[
N_{ij}^k=J_i^k(\pa_pJ^k_j-\pa_jJ^k_p)-J_j^k(\pa_pJ^k_i-\pa_iJ^k_p).
\]
We write $J_i:=J_i^j\pa_j$ for convenience.

\[
N_{ik}^rN_{ri}^sJ_s^k
\]
\[
=\{J_i^p(\pa_pJ^r_k-\pa_kJ^r_p)-J_k^p(\pa_pJ^r_i-\pa_iJ^r_p)\}
\{J_r^q(\pa_qJ^s_i-\pa_iJ^s_q)-J_i^q(\pa_qJ^s_r-\pa_rJ^s_q)\}J_s^k
\]
\[
=\{J_s^kJ_i^p(\pa_pJ^r_k-\pa_kJ^r_p)-J_s^kJ_k^p(\pa_pJ^r_i-\pa_iJ^r_p)\}
\{J_r^q(\pa_qJ^s_i-\pa_iJ^s_q)-J_i^q(\pa_qJ^s_r-\pa_rJ^s_q)\}
\]
\[
=\{(J_s^kJ_i^p\pa_pJ^r_k-J_s^kJ_i^p\pa_kJ^r_p)-J_s^kJ_k^p(\pa_pJ^r_i-\pa_iJ^r_p)\}
\{J_r^q(\pa_qJ^s_i-\pa_iJ^s_q)-J_i^q(\pa_qJ^s_r-\pa_rJ^s_q)\}
\]
\[
=\{J_s^k\cdot J_iJ^r_k-J_i^p\cdot J_sJ^r_p+\pa_sJ^r_i-\pa_iJ^r_s\}
\{J_rJ^s_i-J_r^q\pa_iJ^s_q-J_iJ^s_r+J_i^q\pa_rJ^s_q\}
\]
\[
=\{J_s^k\cdot J_iJ^r_k -\pa_iJ^r_s-J_i^p\cdot J_sJ^r_p +\cdot\pa_sJ^r_i\}
\{-J_iJ^s_r-J_r^q\pa_iJ^s_q+J_rJ^s_i+J_i^q\pa_rJ^s_q\}
\]

\[
=-J_s^k\cdot J_iJ^r_k\cdot J_iJ^s_r-J_r^qJ_s^k\cdot J_iJ^r_k\cdot \pa_iJ^s_q+J_s^k\cdot J_iJ^r_k\cdot J_rJ^s_i+J_i^qJ_s^k\cdot J_iJ^r_k\cdot \pa_rJ^s_q
\]
\[
+J_iJ^s_r\cdot\pa_iJ^r_s+J_r^q\pa_iJ^s_q\cdot\pa_iJ^r_s-J_rJ^s_i\cdot\pa_iJ^r_s-J_i^q\pa_rJ^s_q\cdot\pa_iJ^r_s
\]
\[
+J_i^p\cdot J_sJ^r_p\cdot J_iJ^s_r+J_r^qJ_i^p\cdot J_sJ^r_p\cdot \pa_iJ^s_q-J_i^p\cdot J_sJ^r_p\cdot J_rJ^s_i-J_i^qJ_i^p\cdot J_sJ^r_p\cdot \pa_rJ^s_q
\]
\[
-J_iJ^s_r\cdot\pa_sJ^r_i-J_r^q\pa_iJ^s_q\cdot\pa_sJ^r_i+J_rJ^s_i\cdot\pa_sJ^r_i+J_i^q\pa_rJ^s_q\cdot\pa_sJ^r_i
\]

\[
=-J_s^k\cdot J_iJ^r_k\cdot J_iJ^s_r+J_s^k\cdot J_iJ^r_k\cdot J_rJ^s_i+J_i^p\cdot J_sJ^r_p\cdot J_iJ^s_r-J_i^p\cdot J_sJ^r_p\cdot J_rJ^s_i
\]
\[
-J_r^qJ_s^k\cdot J_iJ^r_k\cdot \pa_iJ^s_q+J_i^qJ_s^k\cdot J_iJ^r_k\cdot \pa_rJ^s_q-J_rJ^s_i\cdot\pa_iJ^r_s
+J_iJ^s_r\cdot\pa_iJ^r_s+J_r^qJ_i^p\cdot J_sJ^r_p\cdot \pa_iJ^s_q
\]
\[
-J_i^qJ_i^p\cdot J_sJ^r_p\cdot \pa_rJ^s_q-J_iJ^s_r\cdot\pa_sJ^r_i+J_rJ^s_i\cdot\pa_sJ^r_i
\]
\[
+J_r^q\pa_iJ^s_q\cdot\pa_iJ^r_s-J_i^q\pa_rJ^s_q\cdot\pa_iJ^r_s
-J_r^q\pa_iJ^s_q\cdot\pa_sJ^r_i+J_i^q\pa_rJ^s_q\cdot\pa_sJ^r_i.
\]

IV3, meaning the third term of fourth line in the above four lines after $=$
\[
-J_r^q\pa_iJ^s_q\cdot\pa_sJ^r_i
=-\pa_sJ^r_i\cdot \pa_iJ^s_q\cdot J_r^q
=J^s_q\cdot \pa_sJ^r_i\cdot \pa_iJ_r^q
=J_qJ^r_i\cdot \pa_iJ_r^q
=J_jJ^k_i\cdot \pa_iJ_k^j.
\]
II3,  meaning the third term of fourth line in the four lines after $=$. Others  are similar.
\[
-J_rJ^s_i\cdot\pa_iJ^r_s=-J_jJ^k_i\cdot\pa_iJ^j_k
\]
\[
\boxed{II3+IV3=0.}
\]

\[
N_{ik}^rN_{ri}^sJ_s^k
\]
\[
=-J_s^k\cdot J_iJ^r_k\cdot J_iJ^s_r+J_s^k\cdot J_iJ^r_k\cdot J_rJ^s_i+J_i^p\cdot J_sJ^r_p\cdot J_iJ^s_r-J_i^p\cdot J_sJ^r_p\cdot J_rJ^s_i
\]
\[
-J_r^qJ_s^k\cdot J_iJ^r_k\cdot \pa_iJ^s_q+J_i^qJ_s^k\cdot J_iJ^r_k\cdot \pa_rJ^s_q+II3
+J_iJ^s_r\cdot\pa_iJ^r_s+J_r^qJ_i^p\cdot J_sJ^r_p\cdot \pa_iJ^s_q
\]
\[
-J_i^qJ_i^p\cdot J_sJ^r_p\cdot \pa_rJ^s_q-J_iJ^s_r\cdot\pa_sJ^r_i+J_rJ^s_i\cdot\pa_sJ^r_i
\]
\[
+J_r^q\pa_iJ^s_q\cdot\pa_iJ^r_s-J_i^q\pa_rJ^s_q\cdot\pa_iJ^r_s
+IV3+J_i^q\pa_rJ^s_q\cdot\pa_sJ^r_i
\]

IV2
\[
-J_i^q\pa_rJ^s_q\cdot\pa_iJ^r_s
=J^s_q\cdot \pa_rJ_i^q\cdot\pa_iJ^r_s
=J^s_q\cdot \pa_jJ_i^q\cdot\pa_iJ^j_s
=-J^j_s\cdot \pa_jJ_i^q\cdot  \pa_i J^s_q
\]
\[
=-J_sJ_i^q\cdot  \pa_i J^s_q=-J_jJ_i^k\cdot  \pa_i J^j_k.
\]

II5
\[
J_r^qJ_i^p\cdot J_sJ^r_p\cdot \pa_iJ^s_q
=-J_i^pJ^s_q\cdot J_sJ^r_p\cdot \pa_iJ_r^q
\]
\[
=-J_i^pJ^s_q\cdot J_s^k\pa_kJ^r_p\cdot \pa_iJ_r^q
=J_i^p\pa_jJ^r_p\cdot \pa_iJ_r^j
=-J^l_k\pa_jJ_i^k\cdot \pa_iJ_l^j
\]
\[
=J_l^j\pa_jJ_i^k\cdot \pa_iJ^l_k=J_lJ_i^k\cdot \pa_iJ^l_k=J_jJ_i^k\cdot \pa_iJ^j_k.
\]

\[
N_{ik}^rN_{ri}^sJ_s^k
\]
\[
=-J_s^k\cdot J_iJ^r_k\cdot J_iJ^s_r+J_s^k\cdot J_iJ^r_k\cdot J_rJ^s_i+J_i^p\cdot J_sJ^r_p\cdot J_iJ^s_r-J_i^p\cdot J_sJ^r_p\cdot J_rJ^s_i
\]
\[
-J_r^qJ_s^k\cdot J_iJ^r_k\cdot \pa_iJ^s_q+J_i^qJ_s^k\cdot J_iJ^r_k\cdot \pa_rJ^s_q+II3
+J_iJ^s_r\cdot\pa_iJ^r_s-J^l_k\pa_jJ_i^k\cdot \pa_iJ_l^j
\]
\[
-J_i^qJ_i^p\cdot J_sJ^r_p\cdot \pa_rJ^s_q-J_iJ^s_r\cdot\pa_sJ^r_i+J_rJ^s_i\cdot\pa_sJ^r_i
\]
\[
+J_r^q\pa_iJ^s_q\cdot\pa_iJ^r_s-J_jJ_i^k\cdot  \pa_i J^j_k
+IV3+J_i^q\pa_rJ^s_q\cdot\pa_sJ^r_i.
\]

II2
\[
+J_i^qJ_s^k\cdot J_iJ^r_k\cdot \pa_rJ^s_q=-J^s_qJ_s^k\cdot J_iJ^r_k\cdot \pa_rJ_i^q=J_iJ^r_k\cdot \pa_rJ_i^k=J_iJ^k_j\cdot \pa_kJ_i^j.
\]
III2
\[
-J_iJ^s_r\cdot\pa_sJ^r_i=-J_iJ^k_j\cdot\pa_kJ^j_i.
\]
\[
\boxed{II2+III2=0.}
\]

\[
\boxed{II5+IV2=0.}
\]

\[
N_{ik}^rN_{ri}^sJ_s^k
\]
\[
=-J_s^k\cdot J_iJ^r_k\cdot J_iJ^s_r+J_s^k\cdot J_iJ^r_k\cdot J_rJ^s_i+J_i^p\cdot J_sJ^r_p\cdot J_iJ^s_r-J_i^p\cdot J_sJ^r_p\cdot J_rJ^s_i
\]
\[
-J_r^qJ_s^k\cdot J_iJ^r_k\cdot \pa_iJ^s_q+II2+II3
+J_iJ^s_r\cdot\pa_iJ^r_s+II5
\]
\[
-J_i^qJ_i^p\cdot J_sJ^r_p\cdot \pa_rJ^s_q+III2+J_rJ^s_i\cdot\pa_sJ^r_i
\]
\[
+J_r^q\pa_iJ^s_q\cdot\pa_iJ^r_s+IV2
+IV3+J_i^q\pa_rJ^s_q\cdot\pa_sJ^r_i.
\]

III1
\[
-J_i^qJ_i^p\cdot J_sJ^r_p\cdot \pa_rJ^s_q=-J^s_qJ^r_p\cdot J_sJ_i^p\cdot \pa_rJ_i^q
=-J^s_q J_s^t\pa_tJ_i^p\cdot J^r_p\pa_rJ_i^q
\]
\[
=\pa_qJ_i^p\cdot J_pJ_i^q=J_jJ_i^k\cdot \pa_kJ_i^j.
\]

II1
\[
-J_r^qJ_s^k\cdot J_iJ^r_k\cdot \pa_iJ^s_q
=J^s_qJ_s^k\cdot J_iJ^r_k\cdot \pa_iJ_r^q=- J_iJ^k_j\cdot \pa_iJ_k^j.
\]
II4
\[
+J_iJ^s_r\cdot\pa_iJ^r_s=+J_iJ^k_j\cdot\pa_iJ^j_k.
\]
\[
\boxed{II1+II4=0.}
\]

\[
N_{ik}^rN_{ri}^sJ_s^k
\]
\[
=-J_s^k\cdot J_iJ^r_k\cdot J_iJ^s_r+J_s^k\cdot J_iJ^r_k\cdot J_rJ^s_i+J_i^p\cdot J_sJ^r_p\cdot J_iJ^s_r-J_i^p\cdot J_sJ^r_p\cdot J_rJ^s_i
\]
\[
II1+II2+II3
+II4+II5
\]
\[
J_jJ_i^k\cdot \pa_kJ_i^j+III2+J_rJ^s_i\cdot\pa_sJ^r_i
\]
\[
+J_r^q\pa_iJ^s_q\cdot\pa_iJ^r_s+IV2
+IV3+J_i^q\pa_rJ^s_q\cdot\pa_sJ^r_i.
\]

IV4
\[
+J_i^q\pa_rJ^s_q\cdot\pa_sJ^r_i=-J^r_i\pa_rJ^s_q\cdot\pa_sJ_i^q
=-J_iJ^s_q\cdot\pa_sJ_i^q=-J_iJ^k_j\cdot\pa_kJ_i^j.
\]

III3
\[
+J_rJ^s_i\cdot\pa_sJ^r_i=+J_jJ^k_i\cdot\pa_kJ^j_i
\]
\[
III1=III3=2J_jJ^k_i\cdot\pa_kJ^j_i.
\]

\[
N_{ik}^rN_{ri}^sJ_s^k
\]
\[
=-J_l^j\cdot J_iJ^k_j\cdot J_iJ^l_k+J_l^j\cdot J_iJ^k_j\cdot J_kJ^l_i+J_i^l\cdot J_iJ^k_j \cdot J_kJ^j_l-J_i^l\cdot J_jJ^k_i \cdot 
J_kJ^j_l\]
\[
II1+II2+II3
+II4+II5
\]
\[
J_jJ_i^k\cdot \pa_kJ_i^j+III2+J_jJ_i^k\cdot \pa_kJ_i^j
\]
\[
+J_r^q\pa_iJ^s_q\cdot\pa_iJ^r_s+IV2
+IV3-J_iJ^k_j\cdot\pa_kJ_i^j.
\]

I3
\[
J_i^l\cdot J_iJ^k_j \cdot J_kJ^j_l=J_i^l\cdot J_iJ^k_j \cdot J_k^p\pa_pJ^j_l=J_k^p\cdot J_iJ^k_j \cdot J_i^l\pa_pJ^j_l
\]
\[
=-J_k^pJ^j_l\cdot J_iJ^k_j \cdot \pa_p J_i^l
=-J^j_l\cdot J_iJ^k_j \cdot J_k^p\pa_p J_i^l
=-J^j_l\cdot J_iJ^k_j \cdot J_k J_i^l.
\]
\[
\boxed{I2+I3=0.}
\]

I4
\[
-J_i^l\cdot J_jJ^k_i \cdot 
J_kJ^j_l=-J_i^l\cdot J_jJ^k_i \cdot 
J_k^p\pa_pJ^j_l
=J^j_l\cdot J_jJ^k_i \cdot 
J_k^p\pa_pJ_i^l
\]
\[
=J^j_l\cdot J_jJ^k_i \cdot J_kJ_i^l
=J^j_l\cdot J_j^p\pa_pJ^k_i \cdot J_kJ_i^l
=- J_kJ_i^l\cdot \pa_lJ^k_i
=- J_jJ_i^k\cdot \pa_kJ^j_i.
\]

\[
\boxed{I4+III1=0,\quad III1=III3.}
\]

\[
N_{ik}^rN_{ri}^sJ_s^k
\]
\[
=-J_l^j\cdot J_iJ^k_j\cdot J_iJ^l_k
+I2+I3
+I4
\]
\[
II1+II2+II3
+II4+II5
\]
\[
III1+III2+J_jJ_i^k\cdot \pa_kJ_i^j
\]
\[
+J_r^q\pa_iJ^s_q\cdot\pa_iJ^r_s+IV2
+IV3-J_iJ^k_j\cdot\pa_kJ_i^j.
\]

\[
N_{ik}^rN_{ri}^sJ_s^k
\]
\[
=-J_l^j\cdot J_iJ^k_j\cdot J_iJ^l_k
+J_jJ_i^k\cdot \pa_kJ_i^j
+J_l^j\pa_iJ^k_j\cdot\pa_iJ^l_k
-J_iJ^k_j\cdot\pa_kJ_i^j
\]
\[
=-J_l^j\cdot J_iJ^k_j\cdot J_iJ^l_k
+J_l^j\pa_iJ^k_j\cdot\pa_iJ^l_k
+J_jJ_i^k\cdot \pa_kJ_i^j
-J_iJ^k_j\cdot\pa_kJ_i^j
\]
\[
=J_l^j( -J_iJ^k_j\cdot J_iJ^l_k
+\pa_iJ^k_j\cdot\pa_iJ^l_k)
+(J_jJ_i^k
-J_iJ^k_j)\cdot\pa_kJ_i^j.
\]

\[
N_{ik}^rN_{ri}^sJ_s^k
\]
\[
=-J_l^j\cdot J_iJ^k_j\cdot J_iJ^l_k
+J_l^j\pa_iJ^k_j\cdot\pa_iJ^l_k
+J_jJ_i^k\cdot \pa_kJ_i^j
-J_iJ^k_j\cdot\pa_kJ_i^j
\]
\[
=-J_l^s\cdot J_pJ^r_s\cdot J_pJ^l_r
+J_l^j\pa_iJ^r_j\cdot\pa_iJ^l_r
+J_pJ_s^r\cdot \pa_rJ_s^p
-J_pJ^r_s\cdot\pa_rJ_p^s
\]
\[
=-J_t^s\cdot J_p^i\pa_iJ^r_s\cdot J_p^j\pa_jJ^t_r
+J_l^j\pa_iJ^r_j\cdot\pa_iJ^l_r
+J_p^i\pa_iJ_s^j\cdot \pa_jJ_s^p
-J_p^i\pa_iJ^j_s\cdot\pa_jJ_p^s
\]

\[
=-J_t^k\cdot J_p^i\pa_iJ^l_k\cdot J_p^j\pa_jJ^t_l
+J_l^j\pa_iJ^k_j\cdot\pa_iJ^l_k
+J_p^i\pa_iJ_k^j\cdot \pa_jJ_k^p
-J_l^i\pa_iJ^j_k\cdot\pa_jJ_l^k
\]
\[
=-J_t^kJ_p^iJ_p^j\cdot \pa_iJ^l_k\cdot \pa_jJ^t_l
+J_l^j\cdot\pa_iJ^k_j\cdot\pa_iJ^l_k
+J_l^i\cdot \pa_iJ_k^j\cdot \pa_jJ_k^l
-J_l^i\cdot\pa_iJ^j_k\cdot\pa_jJ_l^k.
\]

\[
N_{ik}^rN_{ri}^sJ_s^k
\]
\[
=-J_t^kJ_p^iJ_p^j\cdot \pa_iJ^l_k\cdot \pa_jJ^t_l
+J_l^j\cdot\pa_iJ^k_j\cdot\pa_iJ^l_k
+J_l^i\cdot \pa_iJ_k^j\cdot \pa_jJ_k^l
-J_l^i\cdot\pa_iJ^j_k\cdot\pa_jJ_l^k
\]

The first term above is zero.
In fact, 
\[
-J_t^kJ_p^iJ_p^j\cdot \pa_iJ^l_k\cdot \pa_jJ^t_l
\]
\[
=-\pa_i(J_t^kJ_p^iJ_p^j\cdot J^l_k\cdot \pa_jJ^t_l)
\]
\[
+\pa_iJ_t^kJ_p^iJ_p^j\cdot J^l_k\cdot \pa_jJ^t_l
\]
\[
+J_t^k\pa_i J_p^iJ_p^j\cdot J^l_k\cdot \pa_jJ^t_l
\]
\[
+J_t^kJ_p^i\pa_iJ_p^j\cdot J^l_k\cdot \pa_jJ^t_l
\]
\[
+J_t^kJ_p^iJ_p^j\cdot J^l_k\cdot \pa_j\pa_iJ^t_l
\]

\[
=+\pa_i(J_p^iJ_p^j\cdot \pa_jJ^l_l)
\]
\[
+J_p^iJ_p^j\cdot\pa_iJ_t^k\cdot J^l_k\cdot \pa_jJ^t_l
\]
\[
-J_p^j\cdot\pa_i J_p^i\cdot \pa_jJ^l_l
\]
\[
-J_p^i\cdot \pa_iJ_p^j\cdot  \pa_jJ^l_l
\]
\[
-J_p^iJ_p^j\cdot \pa_j\pa_iJ^l_l
\]

\[
=J^l_kJ_p^iJ_p^j\cdot\pa_iJ_t^k\cdot\pa_jJ^t_l
\]

\[
-J_t^kJ_p^iJ_p^j\cdot \pa_iJ^l_k\cdot \pa_jJ^t_l
\]
\[
=J^l_kJ_p^iJ_p^j\cdot\pa_jJ^t_l\cdot\pa_iJ_t^k
\]
\[
=J^b_aJ_p^iJ_p^j\cdot\pa_jJ^c_b\cdot\pa_iJ_c^a
\]
\[
=J^k_tJ_p^iJ_p^j\cdot\pa_jJ^l_k\cdot\pa_iJ_l^t
\]
\[
=J^k_tJ_p^uJ_p^v\cdot\pa_vJ^l_k\cdot\pa_uJ_l^t
\]
\[
=J^k_tJ_p^vJ_p^u\cdot\pa_vJ^l_k\cdot\pa_uJ_l^t
\]
\[
=J^k_tJ_p^iJ_p^j\cdot\pa_iJ^l_k\cdot\pa_jJ_l^t
\]
Therefore
\[
-J_t^kJ_p^iJ_p^j\cdot \pa_iJ^l_k\cdot \pa_jJ^t_l=J^k_tJ_p^iJ_p^j\cdot\pa_iJ^l_k\cdot\pa_jJ_l^t
\]
\[
-J_t^kJ_p^iJ_p^j\cdot \pa_iJ^l_k\cdot \pa_jJ^t_l=0.
\]

Now
\[
N_{ik}^rN_{ri}^sJ_s^k
\]
\[
=J_l^i\cdot \pa_iJ_k^j\cdot \pa_jJ_k^l
-J_l^i\cdot\pa_iJ^j_k\cdot\pa_jJ_l^k
\]
\[
=J_l^i\cdot \pa_iJ_k^j\cdot \pa_j(J_k^l
-J_l^k)
\]
\[
= J_l^i\cdot\pa_j(J_k^l-J_l^k)\cdot\pa_iJ_k^j
\]
\[
=\pa_j (J_l^iJ_k^l-J_l^iJ_l^k)\cdot\pa_iJ_k^j
-\pa_jJ_l^i\cdot(J_k^l-J_l^k)\cdot\pa_iJ_k^j
\]
\[
=\pa_j (-\del^i_k-J_l^iJ_l^k)\cdot\pa_iJ_k^j
-\pa_jJ_l^i\cdot J_k^l\cdot\pa_iJ_k^j+\pa_jJ_l^i\cdot J_l^k\cdot\pa_iJ_k^j
\]

The sum of the last two terms is zero, for
\[
\pa_jJ_l^i\cdot J_l^k\cdot\pa_iJ_k^j
=\pa_iJ_k^j\cdot J_l^k\cdot\pa_jJ_l^i
\]
\[
=\pa_aJ_c^b\cdot J_d^c\cdot\pa_bJ_d^a
=\pa_jJ_l^i\cdot J_k^l\cdot\pa_iJ_k^j.
\]

Therefore
\[
N_{ik}^rN_{ri}^sJ_s^k
=-\pa_j(J_l^iJ_l^k)\cdot\pa_iJ_k^j=L_J
\]
\qed

\end{document}